\documentclass{jnmp}

\usepackage{amsmath}

\JNMPnumberwithin{equation}{section}

\newtheorem{theorem}{Theorem}[section]
\newtheorem{lemma}{Lemma}[section]
\newtheorem{corollary}{Corollary}[section]

\begin{document}

\renewcommand{\evenhead}{D~Roytenberg}
\renewcommand{\oddhead}{Poisson Cohomology of ``Necklace'' Structures}

\thispagestyle{empty}

\FirstPageHead{9}{3}{2002}{\pageref{roytenberg-firstpage}--\pageref{roytenberg-lastpage}}{Article}

\copyrightnote{2002}{D~Roytenberg}

\Name{Poisson Cohomology of $\boldsymbol{SU(2)}$-Covariant
``Necklace'' Poisson Structures on $\boldsymbol{S^{2}}$}
\label{roytenberg-firstpage}

\Author{Dmitry ROYTENBERG}

\Address{Department of Mathematics,
Pennsylvania State University,\\
331 McAllister Building,
University Park, PA 16801,
USA}

\Date{Received December 19, 2001; Accepted April 15, 2002}

\begin{abstract}
\noindent We compute the Poisson cohomology of the one-parameter
family of $ SU(2) $-covariant Poisson structures on the
homogeneous space $ S^{2}={\mathbb C}P^{1}=SU(2)/U(1) $, where $
SU(2) $ is endowed with its standard Poisson--Lie group structure,
thus extending the result of Ginzburg~\cite{Gin1} on the
Bruhat--Poisson structure which is a member of this family. In
particular, we compute several invariants of these structures,
such as the modular class and the Liouville class. As a corollary
of our computation, we deduce that these structures are nontrivial
deformations of each other in the direction of the standard
rotation-invariant symplectic structure on $ S^{2} $; another
corollary is that these structures do not admit smooth rescaling.
\end{abstract}

\section{Introduction}

The Poisson cohomology of a Poisson manifold $ (P,\pi ) $ is the
cohomology of the complex $ (\mathfrak{X} ^{\cdot }(P),d_{\pi
}=[\pi ,\cdot ]) $, where $ \mathfrak{X} ^{k}(P) $ is the space of
smooth $ k $-vector fields on $ P $, and $ [\cdot,\cdot ] $ is the
Schouten bracket. The Poisson cohomology spaces $ H^{k}_{\pi }(P)
$ are important invariants of $ (P,\pi ) $. For instance, $ H_{\pi
}^{0}(P) $ is the space of central (Casimir) functions; $
H^{1}_{\pi }(P) $ is the space of outer derivations of $ \pi $; $
H_{\pi }^{2}(P) $ is the space of non-trivial infinitesimal
deformations of~$ \pi  $, while $ H_{\pi }^{3}(P) $ houses
obstructions to extending a first-order deformation to a formal
deformation. For nondegenerate (symplectic) $ \pi  $, the Poisson
cohomology is isomorphic to the de Rham cohomology of $ P $; in
general, however, this cohomology is notoriously difficult to
compute.

There are two canonical Poisson cohomology classes that merit
special attention. The \emph{modular class $ \Delta \in H^{1}_{\pi
}(P) $} is the obstruction to the existence of an invariant volume
form~\cite{We3}: it vanishes if and only if there exists a measure
on $ P $ preserved by all Hamiltonian flows. The
\emph{Liouville class} is the class of $ \pi  $ itself in $ H_{\pi
}^{2}(P) $. This class is the obstruction to smooth rescaling of $
\pi  $: it vanishes if and only if there exists a vector field~$ X
$ such that $ L_{X}\pi =\pi  $; the flow of this vector field acts
by rescaling $ \pi  $.

The purpose of this note is to compute the Poisson cohomology
 of all $ SU(2) $-covariant
Poisson structures on the two-sphere. Here $ G=SU(2) $ is endowed
with the standard Poisson--Lie group structure and acts on the
homogeneous space $ P=S^{2}=SU(2)/U(1) $ by rotations (recall that
a Lie group $ G $ is a \emph{Poisson--Lie group} if it is endowed
with a~\emph{multiplicative} Poisson tensor, i.e.\ such that the
multiplication $ G\times G\rightarrow G $ is a Poisson map; if a
Poisson--Lie group acts on a manifold $ P $, we say that a Poisson
structure on $ P $ is \emph{$ G $-covariant} if the action $
G\times P\rightarrow P $ is a Poisson map; see \cite{LuWe} or
\cite{ChPr} for details). The $ SU(2) $-covariant structures on $
S^{2} $ form a 1-parameter family $ \pi _{c} $, $ c\in {\mathbb R}
$. For $ |c|>1 $ we get nondegenerate (symplectic) Poisson
structures; $ c=\pm 1 $ corresponds to the (isomorphic) \emph{Bruhat-Poisson}
structures, so called because their symplectic leaves are the Bruhat
cells: a point and an open 2-cell (see \cite{LuWe}); finally, for
$ |c|<1 $ there are two open symplectic leaves (``caps'') of
infinite area separated by a circle (``necklace'') of
zero-dimensional symplectic leaves. All the $ \pi _{c} $'s are
invariant with respect to the residual action of $
S^{1}=U(1)\subset SU(2) $.

The note is organized as follows. Section~\ref{sec:description} is
devoted to the explicit description of the Poisson structures $
\pi _{c} $, while Section~\ref{sec:computation} is devoted to the
computation of their Poisson cohomology, for $ |c|<1 $ (for $
|c|>1 $ it is just the deRham cohomology of $ S^{2} $, whereas the
Bruhat case ($ c=\pm 1 $) was worked out by Viktor
Ginzburg~\cite{Gin1}). We proceed by first linearizing $ \pi _{c}
$ in a neighborhood of the necklace (in an $ S^{1} $-equivariant
way) and computing its local cohomology, then using the
Mayer--Vietoris argument to get the final result, which is
\begin{gather*}
H^{0}_{\pi _{c}}(S^{2})  =  {\mathbb R}, \\ H^{1}_{\pi
_{c}}(S^{2}) =  {\mathbb R}, \\ H^{2}_{\pi _{c}}(S^{2})  =
{\mathbb R} ^{2}
\end{gather*}
 independently of the value of $ c $. In fact, this result coincides with
that of Ginzburg for $ c=1 $. The generator of $ H^{1} $ is the
modular class $ \Delta  $, whereas $ H^{2} $ is spanned by the
classes of $ \pi _{c} $ and $ \pi  $, the inverse of the standard
$ SU(2) $-invariant area form on $ S^{2} $. This shows that (1)~$
\pi _{c} $ does not admit smooth rescaling, and (2)~$ \pi _{c} $
is not isotopic to $ \pi _{c'} $ for $ c\neq c' $.

\section{Description of the Poisson structures\label{sec:description}}

\subsection{The classical $\boldsymbol{r}$-matrix\\
and the standard Poisson--Lie structure on $\boldsymbol{SU(2)}$}

The constructions below can be carried out for any compact
semisimple Lie group, but we will only consider $ SU(2)$.

Recall that the Lie algebra $ \mathfrak {su}(2) $ of $ 2\times 2 $
skew-hermitian traceless matrices has a basis
\[
e_{1}=\frac 12 \left( \begin{array}{cc} i & 0\\ 0 & -i
\end{array}\right) ,\qquad e_{2}=\frac 12 \left( \begin{array}{cc}
0 & 1\\
-1 & 0
\end{array}\right) ,\qquad e_{3}=\frac 12 \left( \begin{array}{cc}
0 & i\\
i & 0
\end{array}\right) \]
 with the commutation relations
 $ [e_{\alpha },e_{\beta }]=\epsilon _{\alpha \beta \gamma }e_{\gamma } $,
where $ \epsilon _{\alpha \beta \gamma } $ is the completely
skew-symmetric symbol. The span of $ e_{1} $ is the Cartan
subalgebra $ \mathfrak {a}\subset \mathfrak {su}(2) $. Recall also
that
\[
SU(2)=\left\{ \left. U=\left( \begin{array}{cc}
u & -\bar{v}\\
v & \bar{u}
\end{array}\right) \right| u,v\in \mathbb C,\ \
 \det U=u\bar{u}+v\bar{v}=1\right\} \]
 identifies $ SU(2) $ with the unit sphere in $ \mathbb C^{2} $.
  The \emph{standard
r-matrix} $ {\textbf {r}}=e_{2}\wedge e_{3}\in \mathfrak
{su}(2)\wedge \mathfrak {su}(2) $ defines a multiplicative Poisson
structure on $ SU(2) $ by
\begin{equation}
\label{PoissonSU(2)} \pi _{SU(2)}(U)={\textbf {r}}U-U{\textbf
{r}}.
\end{equation}
In coordinates,
\begin{gather}
\!\pi \left( \left( \begin{array}{cc} u & -\bar{v}\\ v & \bar{u}
\end{array}\right) \right) =\frac{1}{4}\left( \left( \begin{array}{cc}
v & \bar{u}\\
-u & \bar{v}
\end{array}\right) \wedge \left( \begin{array}{cc}
iv & i\bar{u}\\
iu & -i\bar{v}
\end{array}\right) -\left( \begin{array}{cc}
\bar{v} & u\\
-\bar{u} & v
\end{array}\right) \wedge \left( \begin{array}{cc}
-i\bar{v} & iu\\
i\bar{u} & iv
\end{array}\right) \right)\! \nonumber\\
\qquad{} =-iv\bar{v}\frac{\partial}{\partial u}\wedge
\frac{\partial}{\partial \bar{u}}+\frac 12 \left(
iuv\frac{\partial}{\partial u}\wedge \frac{\partial}{\partial
v}+\overline{iuv\frac{\partial}{\partial u}\wedge
\frac{\partial}{\partial v}}\right) \nonumber\\ \qquad {}+\frac 12
\left( iu\bar{v}\frac{\partial}{\partial u}\wedge
\frac{\partial}{\partial
\bar{v}}+\overline{iu\bar{v}\frac{\partial}{\partial u}\wedge
\frac{\partial}{\partial
\bar{v}}}\right).\label{PoissonSU(2)coord}
\end{gather}
 The Poisson brackets are
\[
\{u,\bar{u}\}=-iv\bar{v}, \qquad \{u,v\}=\frac 12 iuv, \qquad
\{u,\bar{v}\}=\frac 12 iu\bar{v}, \qquad \{v,\bar{v}\}=0.
\]
 It is easy to see that these formulas in fact define a smooth real Poisson
structure on all of $ {\mathbb C} ^{2} $ that restricts to the
unit sphere.

\subsection{The Bruhat--Poisson structure on $\boldsymbol{{\mathbb C} P^{1}}$}

The r-matrix is invariant under the action of the Cartan
subalgebra $ \mathfrak {a} $ , since
\[
[e_{1},{\textbf {r}}]=[e_{1},e_{2}\wedge e_{3}]=[e_{1},e_{2}]
\wedge e_{3}-e_{2}\wedge [e_{1},e_{3}]=e_{3}\wedge
e_{3}+e_{2}\wedge e_{2}=0.
\]
 Hence, the Poisson tensor (\ref{PoissonSU(2)}) vanishes on the maximal torus
(the diagonal subgroup) $ A=U(1)\subset SU(2) $. In particular, $
U(1) $ is a Poisson subgroup, and hence $ \pi _{SU(2)} $ descends
to the quotient $ SU(2)/U(1)=S^{3}/S^{1}=({\mathbb C}
^{2}\setminus 0)/{\mathbb C} ^{\times }={\mathbb C} P^{1}=S^{2} $.
The resulting Poisson structure $ \pi _{1} $ on $ {\mathbb C}
P^{1} $ is called the \emph{Bruhat--Poisson structure} because its
symplectic leaves coincide with the Bruhat cells in $ {\mathbb C}
P^{1} $ \cite{LuWe}: the base point where $ \pi _{1} $ vanishes,
and the complementary open cell where $ \pi _{1} $ is invertible.
It is $ SU(2) $-covariant since $ \pi _{SU(2)} $ is
multiplicative. It is an easy calculation to deduce from
(\ref{PoissonSU(2)coord}) that in the inhomogeneous coordinate
chart $ w=v/u $ covering the base point $ \pi _{1} $ is given by
\[
\pi _{1}=-iw\bar{w}(1+w\bar{w})\frac{\partial}{\partial w}\wedge
\frac{\partial}{\partial\bar{w}}.\]
 In particular, it has a quadratic singularity at $ w=0 $. The
  other inhomogeneous chart $ z=u/v=1/w $ gives coordinates
  on the open symplectic
leaf, in which
\[
\pi _{1}=-i(1+z\bar{z})\frac{\partial}{\partial z}\wedge
\frac{\partial}{\partial \bar{z}}.\]
 The corresponding symplectic 2-form is
\[
\omega _{1}=\frac{idz\wedge d\bar{z}}{1+z\bar{z}}.\] Notice that
this symplectic leaf has infinite area.

\subsection{The other $\boldsymbol{SU(2)}$-covariant Poisson structures
on $\boldsymbol{S^{2}}$}

The difference between any two $ SU(2) $-covariant Poisson
structures on $ {\mathbb C} P^{1} $ is an $ SU(2) $-invariant
bivector field which is Poisson because in two dimensions, any
bivector field is. Thus, any covariant structure is obtained by
adding an invariant structure to the Bruhat structure $ \pi _{1}.
$ To see what these structures look like, it is convenient to
embed the Riemann sphere $ {\mathbb C} P^{1} $as the unit sphere $
S^{2}\subset {\mathbb R} ^{3} $ by the (inverse of) the
stereographic projection. The coordinate transformations are given
by
\begin{gather*}
 x_{1}=\frac{2x}{1+x^{2}+y^{2}}, \qquad x=\frac{x_{1}}{1-x_{3}}, \\
 x_{2}=\frac{2y}{1+x^{2}+y^{2}}, \qquad
y=\frac{x_{2}}{1-x_{3}},\\
x_{3}=\frac{x^{2}+y^{2}-1}{1+x^{2}+y^{2}}, \qquad
x^{2}+y^{2}=\frac{1+x_{3}}{1-x_{3}},
\end{gather*}
 where $ z=x+iy $. We shall identify $ {\mathbb R} ^{3} $
  with $ \mathfrak {su}(2)^{*} $,
with the coadjoint action of $ SU(2) $ by rotations. Then the
linear Poisson structure on $ {\mathbb R} ^{3}=\mathfrak
{su}(2)^{*} $ is given by
\[
-\pi =x_{1}\frac{\partial}{\partial x_{2}}\wedge
\frac{\partial}{\partial x_{3}}+x_{2}\frac{\partial}{\partial
x_{3}}\wedge \frac{\partial}{\partial
x_{1}}+x_{3}\frac{\partial}{\partial x_{1}}\wedge
\frac{\partial}{\partial x_{2}}\]
 whose restriction to the unit sphere (a coadjoint orbit),
  also denoted by $ -\pi  $,
is $ SU(2) $-invariant and symplectic. Moreover, up to a constant
multiple, $ \pi  $ is the only rotation-invariant Poisson
structure on $ S^{2} $: any other invariant structure is of the
form $ \pi '=f\pi  $ for some function $ f $, but since both $ \pi
$ and $ \pi ' $ are invariant, so is $ f $, hence $ f $ is a
constant. It follows that there is a one-parameter family of $
SU(2) $-covariant Poisson structures of the form $ \pi '=\pi
_{1}+\alpha \pi  $, $ \alpha \in {\mathbb R}  $; since $ \pi
_{1}=(1-x_{3})\pi  $ (straightforward calculation), all $ SU(2)
$-covariant structures are of the form
\[
\pi _{c}=\pi _{1}+(c-1)\pi =(c-x_{3})\pi ,\qquad c\in {\mathbb R}.
\]
 It follows that $ \pi _{c} $ is symplectic for $ |c|>1 $, Bruhat for
$ c=\pm 1 $, while for $ |c|<1 $ $ \pi _{c} $ vanishes on the
circle $ \{x_{3}=c\} $ and is nonsingular elsewhere; $ \pi _{c} $
thus has two open symplectic leaves (``caps'') and a ``necklace''
of zero-dimensional symplectic leaves along the circle. It is
these ``necklace'' structures whose Poisson cohomology we shall
compute. Notice that $ \pi _{c} $ and $ \pi _{-c} $ are isomorphic
as Poisson manifolds via $ x_{3}\mapsto -x_{3} $.

In the original $ \{w,\bar{w}\} $-coordinates we have
\begin{gather}
\label{PiStandard} \pi
=-\frac{i}{2}(1+w\bar{w})^{2}\frac{\partial}{\partial w}\wedge
\frac{\partial}{\partial
\bar{w}}=\frac{1}{4}\left(1+x^{2}+y^{2}\right)^{2}\frac{\partial}{\partial
x}\wedge \frac{\partial}{\partial y},\\
 \pi _{c}  =\pi
_{1}+(c-1)\pi =
-\frac{i}{2}(1+w\bar{w})((c+1)w\bar{w}+c-1)\frac{\partial}{\partial
w}\wedge \frac{\partial}{\partial \bar{w}}\nonumber \\
 \phantom{\pi _{c}  =\pi
_{1}+(c-1)\pi}{}
=\frac{1}{4}\left(1+x^{2}+y^{2}\right)\left((c+1)\left(x^{2}+y^{2}\right)
+c-1\right)\frac{\partial}{\partial x}\wedge
\frac{\partial}{\partial y},\label{PiC}
\end{gather}
 where $ w=x+iy $.

\subsection{Symplectic areas and modular vector fields}

Before we proceed to cohomology computations, we shall compute
some invariants of the structures $ \pi _{c} $. For $ |c|>1 $ $
\pi _{c} $ is symplectic, and the only invariant is the symplectic
area. For the other values of $ c $, the areas of the open
symplectic leaves are easily seen to be infinite; instead, we will
compute the modular vector field of $ \pi _{c} $ with respect to
the standard rotation-invariant volume form $ \omega $ on $ S^{2}
$ (the inverse of $ \pi  $). By elementary calculations we obtain
the following

\begin{lemma}
(1) If $ |c|>1, $ the symplectic area of $ (S^{2},\pi _{c}) $ is
given by
\[
V(c)=2\pi \ln \frac{c+1}{c-1}.\]

(2) For all values of $ c $ the modular vector field with respect
to $ \omega  $ is
\[
\Delta _{\omega }=x\frac{\partial}{\partial
y}-y\frac{\partial}{\partial x}.\]
\end{lemma}

\setcounter{corollary}{1}
\begin{corollary}
If $ |c|,\: |c'|>1, $ $ \pi _{c} $ and $ \pi _{c'} $ are not
isomorphic unless $ |c|=|c'| $.
\end{corollary}
\begin{corollary}
\label{Cor:ModClass}If $ |c|<1, $ the modular class of $ \pi _{c}
$ is nonzero.
\end{corollary}

\begin{proof}
The modular vector field $ \Delta _{\omega } $ rotates the
necklace, hence cannot be Hamilto\-nian.~
\end{proof}

In fact, the modular class of the Bruhat--Poisson structures $ \pi
_{\pm 1} $ is also nonzero~\cite{Gin1}.

Unfortunately, the modular vector field does not help us
distinguish the different ``necklace'' structures. The restriction
of $ \Delta _{\omega } $ to the necklace is independent of $
\omega  $ since changing~$ \omega  $ changes $ \Delta _{\omega } $
by a Hamiltonian vector field which necessarily vanishes along the
necklace, so the period of $ \Delta _{\omega } $ restricted to the
necklace is an invariant, but it has the same value of~$ 2\pi $
for all~$ \pi _{c}. $ When we compute the Poisson cohomology of $
\pi _{c} $ we will see a different way to distinguish them.

\section{Computation of Poisson cohomology\label{sec:computation}}

For $ |c|>1 $ $ \pi _{c} $ is symplectic, so its Poisson
cohomology is isomorphic to the de Rham cohomology of $ S^{2} $;
the Poisson cohomology of the Bruhat--Poisson structure $ \pi
_{\pm 1} $ was worked out by Ginzburg~\cite{Gin1}. Here we shall
compute the cohomology of the necklace structures~$ \pi _{c} $ for
$ |c|<1 $. Our strategy will be similar to Ginzburg's: first
compute the cohomology of the formal neighborhood of the necklace,
show that the result is actually valid in a~finite small
neighborhood and finally, use a Mayer--Vietoris argument to deduce
the global result. The validity of the Mayer--Vietoris argument
for Poisson cohomology comes from the simple observation that on
any Poisson manifold $ (P,\pi ) $ the differential $ d_{\pi } $ is
functorial with respect to restrictions to open subsets (i.e.\ a
morphism of the sheaves of smooth multivector fields on $ P $).

It will be convenient to introduce another change of coordinates:
\[
s=\frac{x}{\sqrt{1+x^{2}+y^{2}}}, \qquad
t=\frac{y}{\sqrt{1+x^{2}+y^{2}}}
\]
 mapping the $ (x,y) $-plane to the open unit disk in the $ (s,t) $-plane.
In the new coordinates $ \pi _{c} $ and $ \pi  $ are given by
\begin{gather}
\label{PiC'} \pi _{c}=\frac 12
\left(s^{2}+t^{2}-\frac{1-c}{2}\right)\frac{\partial}{\partial
s}\wedge \frac{\partial}{\partial t},\\ \label{PiStandard'} \pi
=\frac{1}{4}\frac{\partial}{\partial s}\wedge
\frac{\partial}{\partial t}
\end{gather}
 and the necklace is the circle of radius $ R=\sqrt{\frac{1-c}{2}} $. Observe
that rescaling $ s=\alpha s'$, $t=\alpha t' $ ($ \alpha >0 $)
takes $ \pi _{c} $ with necklace radius $ R $ to $ \pi _{c'} $
with necklace radius $ R'=R/\alpha  $. But this is only a local
isomorphism: it does not extend to all of $ S^{2} $ since it is
not a diffeomorphism of the unit disk. In any case, it shows that
all necklace structures are locally isomorphic, so for local
computations we may assume that $ \pi _{c} $ is given in suitable
coordinates by
\[
\pi _{c}=\frac 12
\left(s^{2}+t^{2}-1\right)\frac{\partial}{\partial s}\wedge
\frac{\partial}{\partial t}.\]

\subsection{Cohomology of the formal neighborhood of the necklace}

Since $ \pi _{c} $ is rotation-invariant, we can lift the
computations in the formal neighborhood of the unit circle in the
$ (s,t) $-plane to its universal cover by introducing
``action-angle coordinates'' $ (I,\theta ) $:
\[
s=\sqrt{1+I}\cos \theta, \qquad t=\sqrt{1+I}\sin \theta
\]
 in which $ \pi _{c} $ is linear:
\[
\pi _{c}=I\frac{\partial}{\partial I}\wedge
\frac{\partial}{\partial \theta }.\]
 Of course we will have to restrict attention to
 multivector fields whose coefficients are periodic in $ \theta  $.
 It will be convenient to think
of multivector fields as functions on the supermanifold with
coordinates $ (I,\theta ,\xi ,\eta ) $ where $ \xi  $``=''$
\partial _{I} $ and $ \eta  $``=''$ \partial _{\theta } $ are
Grassmann (anticommuting) variables. Then $ \pi _{c}=I\xi \eta $
is a function and
\[
d_{\pi _{c}}=[\pi _{c},\cdot ] =-I\eta \frac{\partial}{\partial
I}+I\xi \frac{\partial}{\partial \theta }-\xi \eta
\frac{\partial}{\partial \xi }\]
 is a (homological) vector field. Since $ d_{\pi _{c}} $
  commutes with rotations, we can split the complex into Fourier modes
\[
\mathfrak {X}_{n}^{0}=\big\{f(I)e^{in\theta }\big\};\qquad
\mathfrak {X}_{n}^{1}=\big\{(f(I)\xi +g(I)\eta )e^{in\theta
}\big\};\qquad \mathfrak {X}_{n}^{2}=\big\{h(I)\xi \eta
e^{in\theta }\big\},
\]
 where $ f(I) $, $ g(I) $ and $ h(I) $ are formal power series in $ I $.
It will be convenient to treat the zero and non-zero modes separately; it will
turn out that the cohomology is concentrated entirely in the zero mode.

\medskip

\noindent {\it Case 1.} \textbf{The zero mode
($\boldsymbol{n=0}$)} consists of multivector fields independent
of $ \theta  $, so $ d_{\pi _{c}} $ becomes
\[
\left. d_{\pi _{c}}\right| _{\mathfrak {X}_{0}}=-I\eta
\frac{\partial}{\partial I}+\eta \xi \frac{\partial}{\partial \xi
}\]
 which preserves the degree in $ I $ so the complex $ \mathfrak {X}_{0} $ splits
further into a direct product of sub-complexes $ \mathfrak
{X}_{0,m} $, $ m\geq 0 $ according to the degree:
\[
0\rightarrow \mathfrak {X}^{0}_{0,m}\rightarrow \mathfrak
{X}^{1}_{0,m}\rightarrow \mathfrak {X}^{2}_{0,m}\rightarrow 0.\]
 These complexes are very small ($ \mathfrak {X}_{0,m}^{0} $
  and $ \mathfrak {X}_{0,m}^{2} $
are one-dimensional, while $ \mathfrak {X}_{0,m}^{2} $ is
two-dimensional) and their cohomology is easy to compute. For $
f=cI^{m}\in \mathfrak {X}_{0,m}^{0} $, $ d_{\pi
_{c}}f=-cmI^{m}\eta $, while for $ X=aI^{m}\xi +bI^{m}\eta \in
\mathfrak {X}_{0,m}^{1} $, $ d_{\pi _{c}}X=a(m-1)I^{m}\xi \eta  $.
Therefore, it is clear that for $ m>1 $ the complex is acyclic. On
the other hand, the cohomology of $ \mathfrak {X}_{0,0} $ is
generated by $ 1\in \mathfrak {X}_{0,0}^{0} $ and $ \eta \in
\mathfrak {X}_{0,0}^{1} $, while the cohomology of $ \mathfrak
{X}_{0,1} $ is generated by $ I\xi \in \mathfrak {X}_{0,1}^{1} $
and $ I\xi \eta \in \mathfrak {X}_{0,1}^{2} $. Putting these
together we obtain
\begin{gather}
H_{0}^{0}  =  {\mathbb R}   =  \textrm{span}\{1\},\nonumber\\
H_{0}^{1} =  {\mathbb R} ^{2}  =  \textrm{span}\{\partial _{\theta
},I\partial _{I}\},\nonumber\\ H_{0}^{2}  =  {\mathbb R}   =
\textrm{span}\{I\partial _{I}\wedge \partial _{\theta }\}.
\label{InvHom}
\end{gather}

\noindent {\it Case 2.} \textbf{The non-zero modes
$\boldsymbol{(n\neq 0)}$.} In this case $ d_{\pi _{c}} $ does not
preserve the $ I $-grading so we'll have to consider all power
series at once. Let
\begin{gather*}
f  =  \left(\sum ^{\infty }_{m=0}f_{m}I^{m}\right)e^{in\theta }
\in  \mathfrak {X}_{n}^{0},\\ X  = \left(\sum ^{\infty
}_{m=0}a_{m}I^{m}\right)e^{in\theta }\xi +\left(\sum ^{\infty
}_{m=0}b_{m}I^{m}\right)e^{in\theta }\eta   \in   \mathfrak
{X}_{n}^{1},\\ B  = \left (\sum ^{\infty
}_{m=0}c_{m}I^{m}\right)e^{in\theta }\xi \eta   \in   \mathfrak
{X}_{n}^{2}.
\end{gather*}
 Then
\begin{gather*}
d_{\pi _{c}}f  = \left(\sum ^{\infty
}_{m=1}inf_{m-1}I^{m}\right)e^{in\theta }\xi +\left(\sum ^{\infty
}_{m=1}mf_{m}I^{m}\right)e^{in\theta }\eta, \\ d_{\pi _{c}}X  =
\left(-a_{0}+\sum ^{\infty
}_{m=1}((m-1)a_{m}+inb_{m-1})I^{m}\right)e^{in\theta }\xi \eta
\end{gather*}
 (and, of course, $ d_{\pi _{c}}B=0 $). We see immediately
 that $ d_{\pi _{c}}f=0\Leftrightarrow f=0 $,
hence $ H_{n}^{0}=\{0\} $. Moreover, any $ B $ is a coboundary:
\[
B=d_{\pi _{c}}\left( \left(\sum ^{\infty }_{m\neq
1}\frac{c_{m}}{m-1}I^{m}\right)e^{in\theta }\xi
+\frac{c_{1}}{in}e^{in\theta }\eta \right) \]
 so $ H_{n}^{2}=\{0\} $ as well. Now, $ X $ is a cocycle if and only if
\[
a_{0}  =  b_{0}  =0, \qquad b_{m}  =  -\frac{ma_{m+1}}{in}, \quad
m\geq 1.\]
 Let $ f_{m}=\frac{a_{m+1}}{in} $ for $ m\geq 0 $, $ f=\sum f_{m}I^{m} $.
Then $ X=d_{\pi _{c}}f $. Hence $ H_{n}^{1} $ is also trivial. So
for $ n\neq 0 $ $ \mathfrak {X}_{n} $ is acyclic.

\medskip

It follows that the Poisson cohomology of the formal neighborhood
of the necklace is as in (\ref{InvHom}).

\subsection{Justification for the smooth case}

To see that the cohomology of a finite small neighborhood of the
necklace is the same as for the formal neighborhood we apply an
argument similar to Ginzburg's~\cite{Gin1}. For each Fourier mode
consider the following exact sequence of complexes:
\[
0\rightarrow \mathfrak {X}^{\star }_{n,\textrm{flat }}\rightarrow
\mathfrak {X}^{\star }_{n,\textrm{smooth }}\rightarrow \mathfrak
{X}^{\star }_{n,\textrm{formal }}\rightarrow 0,\]
 where $ \mathfrak {X}^{\star }_{n,\textrm{flat }} $
 consists of smooth multivector
fields whose coefficients vanish along the necklace together with
all derivatives. This sequence is exact by a theorem of E~Borel.
It suffices to show that the flat complex is acyclic. But $ \pi
_{c}^{\#}:\mathfrak {X}^{\star }_{n,\textrm{flat }}\rightarrow
\Omega ^{\star }_{n,\textrm{flat}} $ is an isomorphism since the
coefficient of $ \pi _{c} $ is a polynomial in $ I $, and every
flat form can be divided by a polynomial with a flat result.
Furthermore, the flat deRham complex is acyclic by the homotopy
invariance of deRham cohomology.

Finally, we observe that a smooth multivector field in a
neighborhood of the necklace (given by a \emph{convergent} Fourier
series) is a coboundary if and only if each mode is, and the
primitives can be chosen so that the resulting series converges,
as can be seen from the calculations in the previous subsection
(integration can only improve convergence). Therefore, the Poisson
cohomology of an annular neighborhood $ U $ of the necklace is
\begin{gather}
 H_{\pi _{c}}^{0}(U)  =  {\mathbb R}   =
\textrm{span}\{1\},\nonumber\\ H_{\pi _{c}}^{1}(U)  =  {\mathbb R}
^{2}  =  \textrm{span}\{\partial _{\theta },I\partial
_{I}\},\nonumber\\ H_{\pi _{c}}^{2}(U)  =  {\mathbb R}   =
\textrm{span}\{I\partial _{I}\wedge
\partial _{\theta }\}.
\label{LocHom}
\end{gather}
 Notice that the generators of $ H_{\pi _{c}}^{1}(U) $
 are the rotation $ \partial _{\theta }=s\partial _{t}-t\partial _{s} $
(the modular vector field) and the dilation $ I\partial
_{I}=\frac{s^{2}+t^{2}-1}{2(s^{2}+t^{2})}(s\partial _{s}+t\partial
_{t}) $, while the generator of $ H^{2}_{\pi _{c}}(U) $ is $ \pi
_{c} $ itself, so in particular $ \pi _{c} $ does not admit
rescalings even locally.

\subsection{From local to global cohomology}

We now have all we need to compute the Poisson cohomology of a
necklace Poisson structure $ \pi _{c} $ on $ S^{2} $. Cover $
S^{2} $ by two open sets $ U $ and $ V $ where $ U $ is an annular
neighborhood of the necklace as above, and $ V $ is the complement
of the necklace consisting of two disjoint open caps, on each of
which $ \pi _{c} $ is nonsingular, so that the Poisson cohomology
of $ V $ and $ U\cap V $ is isomorphic to the deRham cohomology.
The short exact Mayer--Vietoris sequence associated to this cover
\[
0\rightarrow \mathfrak {X}^{\star }(S^{2})\rightarrow \mathfrak
{X}^{\star }(U)\oplus \mathfrak {X}^{\star }(V)\rightarrow
\mathfrak {X}^{\star }(U\cap V)\rightarrow 0\]
 leads to a long exact sequence in cohomology:
\[
\begin{array}{ccccccccc}
0 & \rightarrow  & H^{0}_{\pi _{c}}(S^{2}) & \rightarrow  & H_{\pi _{c}}^{0}(U)\oplus H^{0}_{\pi _{c}}(V) & \rightarrow  & H^{0}_{\pi _{c}}(U\cap V) & \rightarrow  & \\
 & \rightarrow  & H_{\pi _{c}}^{1}(S^{2}) & \rightarrow  & H_{\pi _{c}}^{1}(U)\oplus H^{1}_{\pi _{c}}(V) & \rightarrow  & H^{1}_{\pi _{c}}(U\cap V) & \rightarrow  & \\
 & \rightarrow  & H_{\pi _{c}}^{2}(S^{2}) & \rightarrow
 & H_{\pi _{c}}^{2}(U)\oplus H^{2}_{\pi _{c}}(V) & \rightarrow
  & H^{2}_{\pi _{c}}(U\cap V) & \rightarrow  & 0.
\end{array}\]
 Now, the first row is clearly exact since a Casimir function on $ S^{2} $
must be constant on each of the two open symplectic leaves
comprising $ V $, hence constant on all of $ S^{2} $ by
continuity. On the other hand, $ H^{1}_{\pi _{c}}(V)=H^{2}_{\pi
_{c}}(V)=H^{2}_{\pi _{c}}(U\cap V)=\{0\} $. Combining this with
(\ref{LocHom}), we see that what we have left is
\[
\begin{array}{cccccccc}
 &  &  &  & {\mathbb R} ^{2} &  & {\mathbb R} ^{2} & \\
 &  &  &  & \Vert  &  & \Vert  & \\
0 & \rightarrow  & H_{\pi _{c}}^{1}(S^{2}) & \rightarrow  & H_{\pi _{c}}^{1}(U)\oplus H^{1}_{\pi _{c}}(V) & \rightarrow  & H^{1}_{\pi _{c}}(U\cap V) & \rightarrow \\
 &  &  &  &  &  &  & \\
 & \rightarrow  & H_{\pi _{c}}^{2}(S^{2}) & \rightarrow  & H_{\pi _{c}}^{2}(U)\oplus H^{2}_{\pi _{c}}(V) & \rightarrow  & 0. & \\
 &  &  &  & \Vert  &  &  & \\
 &  &  &  & {\mathbb R}  &  &  &
\end{array}\]
 Now, on the one hand, we know by Corollary \ref{Cor:ModClass}
 that $ H_{\pi _{c}}^{1}(S^{2}) $
is at least one-dimensional; on the other hand, the restriction of
the dilation vector field $ I\partial _{I} $ to $ U\cap V $ is not
Hamiltonian: it corresponds under $ \pi _{c}^{\#} $ to the
generator of the first de Rham cohomology of the annulus, diagonally
embedded into $ H^1(U\cap V) $ (a disjoint union of two annuli). It
follows that $ H_{\pi _{c}}^{1}(S^{2}) $ is exactly
one-dimensional, while $ H_{\pi _{c}}^{2}(S^{2}) $ is
two-dimensional.

It only remains to identify the generators. $ H_{\pi
_{c}}^{1}(S^{2}) $ is generated by the modular class, while one of
the generators of $ H_{\pi _{c}}^{2}(S^{2}) $ is $ \pi _{c} $
itself, since its class was shown to be nontrivial even locally.
The other generator is the image of $ (I\partial _{I},-I\partial
_{I})\in H_{\pi _{c}}^{1}(U\cap V) $ under the connecting
homomorphism. This is somewhat unwieldy since it involves a
partition of unity subordinate to the cover $ \{U,V\} $ which does
not yield a clear geometric interpretation of the generator.
Instead, we will show directly that the standard rotationally
invariant symplectic Poisson structure $ \pi  $ on $ S^{2} $ is
nontrivial in $ H_{\pi _{c}}^{2}(S^{2}) $ and so can be taken as
the second generator.

\begin{lemma}
The class of the standard $ SU(2) $-invariant Poisson structure $
\pi  $ on $ S^{2} $ is nonzero in $ H_{\pi _{c}}^{2}(S^{2}) $.
\end{lemma}

\begin{proof}
We will work in coordinates $ (s,t) $ on the unit disk in which $
\pi  $ and $ \pi _{c} $ are given, respectively by
(\ref{PiStandard'}) and(\ref{PiC'}). Locally $ \pi  $ is a
coboundary whose primitive is given by an Euler vector field $
E=\frac{1}{2(c-1)}(s\partial _{s}+t\partial _{t}) $: it's easy to
check that $ [\pi _{c},E]=\pi  $. But $ E $ does not extend to a
vector field on $ S^{2} $ since it does not behave well ``at
infinity'', i.e.\ on the unit circle in the $ (s,t) $-plane.
Therefore, to prove that $ \pi  $ is globally nontrivial it
suffices to show that there does not exist a Poisson vector field
$ X $ such that $ E+X $ is tangent to the unit circle and the
restriction is rotationally invariant. In fact, it suffices to
show that there is no Hamiltonian vector field $ X_{f} $ such that
$ E+X_{f} $ vanishes on the unit circle (since we can always add a
multiple of the modular vector field to cancel the rotation).
Assuming that such an $ f $ exists, we will have, in the polar
coordinates $ s=r\cos \phi  $, $ t=r\sin \phi  $:
\[
E+X_{f}=\frac{1}{2(c-1)}r\frac{\partial}{\partial r}+\frac{1}{2r}
\left(r^{2}-\frac{1-c}{2}\right)\left( \frac{\partial f}{\partial
\phi }\frac{\partial}{\partial r}-\frac{\partial f}{\partial
r}\frac{\partial}{\partial \phi }\right).
\]
 Upon restriction to $ r=1 $ this becomes
\[
\left. \left( E+X_{f}\right) \right| _{r=1} =\left(
\frac{1}{2(c-1)}+\frac{c+1}{4}\left. \frac{\partial f}{\partial
\phi }\right| _{r=1}\right) \left. \frac{\partial}{\partial
r}\right| _{r=1}+\frac{c+1}{4}\left. \frac{\partial f}{\partial
r}\right| _{r=1}\left. \frac{\partial}{\partial \phi }\right|
_{r=1}.\]
 In order for this to vanish it is necessary, in particular,
 that $ \left. \frac{\partial f}{\partial \phi }\right| _{r=1} $
be a nonzero constant which is impossible since $ f $ is periodic
in $ \phi .$
\end{proof}
We have now arrived at our final result:

\setcounter{theorem}{1}

\begin{theorem}
The Poisson cohomology of a necklace Poisson structure $ \pi _{c}
$ on $ S^{2} $ is given as follows:
\begin{gather*}
H_{\pi _{c}}^{0}(S^{2})  =  {\mathbb R}   = {\rm span}\{1\},\\
H_{\pi _{c}}^{1}(S^{2})  =  {\mathbb R}   = {\rm span}\{\Delta
_{\omega }\},\\ H_{\pi _{c}}^{2}(S^{2})  =  {\mathbb R} ^{2}  =
{\rm span}\{\pi _{c},\pi \}.
\end{gather*}
\end{theorem}

\setcounter{corollary}{2}

\begin{corollary}
$ \pi _{c} $ does not admit smooth rescaling.
\end{corollary}

\begin{corollary}
The necklace structures $ \pi _{c} $ and $ \pi _{c'} $ for $ c\neq
c' $ are nontrivial deformations of each other.
\end{corollary}

\begin{proof}
$ \pi _{c'}-\pi _{c} $ is a nonzero multiple of $ \pi  $ but $ \pi
$ is nontrivial in $ H_{\pi _{c}}^{2}(S^{2}) $.
\end{proof}

\subsection*{Acknowledgements}
This work was carried out in the Spring of 1998 at UC Berkeley as
part of the author's dissertation research, and became a part of
his Ph.D.~thesis~\cite{Roy1}. The author wishes to thank his
advisor, Professor Alan Weinstein, for his generous help,
encouragement and support. The author's gratitude also goes to the
Alfred P~Sloan Foundation for financial support throughout this
project.

\label{roytenberg-lastpage}


\begin{thebibliography}{99}\small

\bibitem{ChPr}
Chari V and Pressley~A, A Guide to Quantum Groups, Cambridge Univ.
Press, 1994.

\bibitem{Gin1}
Ginzburg V~L, Momentum Mappings and Poisson Cohomology, {\it Int.
J. Math.}, {\bf 7}, Nr.~3 (1986), 329--358.

\bibitem{LuWe}
 Lu J-H and Weinstein~A,
Poisson {Lie} Groups, Dressing Transformations and Bruhat
  Decompositions, {\it J. Diff. Geom.} {\bf 31} (1990), 501--526.

\bibitem{Roy1}
Roytenberg D, Courant Algebroids, Derived Brackets and Even
Symplectic Supermanifolds, PhD thesis, UC Berkeley, 1999
[math.DG/9910078].

\bibitem{We3}
Weinstein~A, The Modular Automorphism Group of a Poisson Manifold,
{\it J. Geom. and Phys.} {\bf 23} (1997), 379--394.

\end{thebibliography}
\end{document}